\documentclass[a4paper,11pt,leqno,twoside]{article}
\usepackage{amsmath,amsthm,amssymb}
\numberwithin{equation}{section}

\newcommand{\Area}{\mathrm{Area\,}}
\newcommand{\pp}{\sqrt{+^2}}
\newcommand{\mb}{\mathbf}
\newcommand{\bsym}{\boldsymbol}
\newcommand{\bslO}{\backslash\{\mb{0}\}}
\newcommand{\w}[1]{\widetilde{#1}}

\makeatletter

\ifcase\@ptsize
\def\rpkern{\mathchoice{\kern-1.45em}{\kern-1.11em}{}{}}%
\def\grpkern{\mathchoice{\kern-1.013em}{\kern-0.825em}{}{}}%
\or
\def\rpkern{\mathchoice{\kern-1.44em}{\kern-1.11em}{}{}}%
\def\grpkern{\mathchoice{\kern-1.00em}{\kern-0.81em}{}{}}%
\or
\def\rpkern{\mathchoice{\kern-1.472em}{\kern-1.14em}{}{}}%
\def\grpkern{\mathchoice{\kern-1.00em}{\kern-0.815em}{}{}}%
\fi

\def\minibullet{\mathchoice%
{\raise0.2ex\hbox{$\scriptstyle\bullet$}}%
{\raise0.26ex\hbox{$\scriptscriptstyle\bullet$}}{}{}}
\def\butabullet{\mathchoice%
{\raise0.8ex\hbox{$\scriptstyle\bullet$}{\kern-0.365em}%
\lower0.4ex\hbox{$\scriptstyle\bullet$}}%
{\raise0.75ex\hbox{$\scriptscriptstyle\bullet$}{\kern-0.335em}%
\lower0.25ex\hbox{$\scriptscriptstyle\bullet$}}{}{}}

\def\regprod{\operatorname{\coprod\rpkern\prod}\displaylimits}

\def\customprod#1#2%
{\operatorname{\coprod\kern#2{#1}\kern#2\prod}\displaylimits}

\newcommand{\bC}{\mathbb{C}}

\newcommand{\bN}{\mathbb{N}}

\newcommand{\bQ}{\mathbb{Q}}
\newcommand{\bR}{\mathbb{R}}
\newcommand{\bS}{\mathbb{S}}
\newcommand{\bT}{\mathbb{T}}

\newcommand{\bZ}{\mathbb{Z}}

\newcommand{\cL}{\mathcal{L}}

\renewcommand{\a}{\alpha}
\renewcommand{\b}{\beta}
\newcommand{\g}{\gamma}
\renewcommand{\d}{\delta}
\newcommand{\e}{\varepsilon}
\newcommand{\z}{\zeta}

\renewcommand{\th}{\theta}

\renewcommand{\l}{\lambda}
\newcommand{\m}{\mu}
\newcommand{\n}{\nu}
\newcommand{\x}{\xi}
\renewcommand{\r}{\rho}

\renewcommand{\t}{\tau}
\renewcommand{\c}{\chi}

\newcommand{\G}{\Gamma}
\newcommand{\D}{\Delta}

\newcommand{\Om}{\Omega}

\newcommand{\Thm}[1]{{\textrm{Theorem}\,\ref{thm:#1}}}

\newcommand{\Lem}[1]{{\textrm{Lemma}\,\ref{lem:#1}}}
\newcommand{\Cor}[1]{{\textrm{Corollary}\,\ref{cor:#1}}}

\renewcommand{\Re}{\mathrm{Re}\,}

\renewcommand{\det}{\mathrm{det}\,}

\newcommand{\Prim}{\mathrm{Prim}\,}

\newcommand{\bslo}{\backslash\{0\}}
\newcommand{\p}{\partial}

\newcommand{\boldtitle}[1]{\title{\bfseries #1}}

\newcommand{\layout}{%
\topmargin -1in
\oddsidemargin -1in
\evensidemargin -1in
\hoffset -0.07in
\voffset 0.00in
\newlength{\sidemargin}
\setlength{\sidemargin}{-2in}
\addtolength{\sidemargin}{\paperwidth}
\addtolength{\sidemargin}{-\textwidth}
\divide \sidemargin by 2
\setlength{\oddsidemargin}{\sidemargin}
\setlength{\evensidemargin}{\sidemargin}
\setlength{\topmargin}{-2in}
\addtolength{\topmargin}{\paperheight}
\addtolength{\topmargin}{-\textheight}
\addtolength{\topmargin}{-\headheight}
\addtolength{\topmargin}{-\headsep}
\divide \topmargin by 2}

\newenvironment{MSC}{%
\smallbreak
\noindent 2000\ \textbf{Mathematics Subject Classification\,:}\ 
{\itshape Primary}}

\newenvironment{keywords}{%
\noindent\textbf{Key words and phrases\,:}\itshape}

\theoremstyle{theorem}
\newtheorem*{multitheorem}{\variable@name}

\theoremstyle{plain}
\newtheorem{thm}{Theorem}[section]
\newtheorem{prop}[thm]{Proposition}
\newtheorem{lem}[thm]{Lemma}
\newtheorem{cor}[thm]{Corollary}

\newtheorem{conj}[thm]{Conjecture}

\theoremstyle{definition}

\newtheorem{example}[thm]{Example}
\newtheorem{remark}[thm]{Remark}

\boldtitle{Ruelle type $L$-functions versus determinants of Laplacians
for torsion free abelian groups\thanks{Partially supported by
Grant-in-Aid for Scientific Research (B) No.15340012}}
\author{Nobushige KUROKAWA, Masato WAKAYAMA and\\ Yoshinori YAMASAKI}
\date{\today}

\layout

\begin{document}

\setlength{\baselineskip}{15pt}
\maketitle

\begin{abstract}
 We study Ruelle's type zeta and $L$-functions for a torsion free
 abelian group $\G$ of rank $\n\ge 2$ defined via an Euler product. It
 is shown that the imaginary axis is a natural boundary of this zeta
 function when $\n=2,4$ and $8$, and in particular, such a zeta function
 has no determinant expression. Thus, conversely, expressions like
 Euler's product for the determinant of the Laplacians of the torus
 $\bR^{\n}/\G$ defined via zeta regularizations are investigated. Also,
 the limit behavior of an arithmetic function arising from the Ruelle
 type zeta function is observed.  
\begin{MSC}
 11M36, 11N37
\end{MSC} 
\end{abstract}

\section{Introduction}

 The Ruelle zeta function is a particular kind of dynamical zeta
 function which counts periodic orbits for maps or flows. It was
 introduced by Ruelle in the mid-1970s \cite{Ruelle1976}. The Ruelle
 zeta function is also regarded (essentially) as the Selberg zeta
 function when one considers the primitive geodesic flows on a
 hyperbolic space. 

 Instead, we have introduced in \cite{TorusSelbergZeta} various Ruelle
 type zeta functions for a complex torus and studied asymptotic averages
 of certain arithmetic functions arising from the logarithm of such
 Ruelle type zeta functions. In this paper, we study analytic properties
 of a Ruelle type $L$-function for a general torsion free abelian group
 and calculate also the determinant of the Laplacian for the
 corresponding group.

 Let $\G$ be a torsion free abelian group of rank $\n$ $(\n\ge 2)$. Let
 $\t_1,\ldots,\t_{\n}$ be a generator of $\G$. Then we may write 
 $\G=\bZ\t_1\oplus\cdots\oplus\bZ\t_{\n}$. A non-zero element 
 $P=\sum^{\n}_{j=1} n_j\t_j\in \G$ is said to be primitive if the
 integers $n_j$ are relatively prime, that is, the greatest common
 divisor $\gcd(n_1,\ldots,n_{\n})$ is equal to $1$. Here we interpret
 $\gcd(n_1,\ldots,n_{\n})$ as the greatest common divisor of the all
 positive factors of $|n_j|$ $(j=1,2,\ldots,\n)$. We denote by
 $\Prim(\G)$ the set of all primitive elements in $\G$. Let
 $\ell:\G \to \bR_{\ge 0}$ be a function satisfying the homogeneity
 condition $\ell(j\g)=j\ell(\g)$ for any $j\in\bN$ and $\g\in \G$. For
 each element $\g\in \G$ we define the norm $N(\g)$ by
 $N(\g)=e^{\ell(\g)}$. Note that any element
 $\g=\sum^{\n}_{j=1} n_j\t_j\in \G$ can be uniquely expressed as 
 $\g=dP$, where $d=\gcd(n_1,\ldots,n_{\n})$ and $P=P_{\g}\in\Prim(\G)$.
 It is clear that $N(\g)=N(P_\g)^d$. Further, let $\r:\G \to U(N)$ be an
 $N$-dimensional unitary representation of $\G$. For a given such norm
 function $N(\g)$ (or the length function $\ell(\g)$) and a
 representation $\r$ of $\G$, we introduce a Ruelle type $L$-function
 $L_\G(s,\r)$ for $\G$ by the Euler product as
\[
 L_\G(s,\r)=L_\G(s,\r;\n)
:=\prod_{P\in \Prim(\G)}\det(1-\r(P)N(P)^{-s})^{-1}.
\]
 When $\r$ is the trivial representation of $\G$, we write
 $\z_\G(s;\n)=L_\G(s,\text{the trivial rep.};\n)$ and call it a Ruelle
 type zeta function for $\G$. Since $\G$ is abelian, any irreducible
 unitary representation of $\G$ is one dimensional. Thus, if we write the
 irreducible decomposition of $\r$ as $\r=\c_1\oplus\cdots\oplus\c_N$
 (here we are taking account the multiplicity in the decomposition)
 where $\c_j$ are one dimensional characters of $\G$, we easily see that  
 $L_\G(s,\r)=\prod_{P\in\Prim(\G)}\prod^N_{j=1}(1-\c_j(P)N(P)^{-s})^{-1}=\prod^{N}_{j=1}L_\G(s,\c_{j})$. 
 In view of this decomposition, it is enough to treat the case where
 $\r$ is one dimensional in our analysis.

 The aim of the present paper is to show first that, for $\n=2,4$ or
 $8$, the Ruelle type zeta function $\z_{\G}(s;\n)$ has a natural
 boundary $\Re(s)=0$ when $\ell(\g)$ is defined by the usual Euclidean 
 norm (\Thm{nb}). This shows, in particular, that such $\z_{\G}(s;\n)$
 does not have a determinant expression for the Laplacian of the torus
 $\bR^{\n}/\G$ (here the determinant is defined by the zeta
 regularization. See \S 3). Thus, conversely, we study the determinant
 $\det{(\D+s^2)}$ of the Laplacian $\D$ on the torus and show that it
 has a sort of Euler products (\Thm{DetLaplacian}). Though the results
 show that $\det{(\D+s^2)}$ can not give $\z_{\G}(s;\n)$, we notice that
 the logarithm of both the functions involve certain common arithmetic
 function (compare \Cor{series} with \eqref{for:logLpp}). Therefore, we
 also discuss asymptotic averages of certain arithmetic functions
 arising from these Ruelle type zeta functions like in
 \cite{TorusSelbergZeta}.

\section{$\bsym{L_{\pp}(s,\c)}$}

 We study a Ruelle type $L$-function $L_\G(s,\r;\n)$ when
 $\ell(\g)=\sqrt{{n_1}^2+\cdots+{n_{\n}}^2}$ for
 $\g=\sum^{\n}_{j=1}n_j\t_j\in\G$. To distinguish the case from a general
 choice of length functions, we write $L_{\pp}(s,\r;\n)$ instead of
 $L_\G(s,\r;\n)$. Let $\c$ be a unitary character of
 $\G=\bZ\t_1\oplus\cdots\oplus\bZ\t_{\n}$. Put
 $\c(\t_j)=e^{2\pi i\a_j}$ for $j=1,\ldots,\n$. Then we have
 $\c(\sum^{\n}_{j=1}n_j\t_j)=e^{2\pi i(n_1\a_1+\cdots+n_{\n}\a_{\n})}$.
 To economize the space, we use a multi-index. Write an element
 $\g=\sum^{\n}_{j=1}n_j\t_j\in\G$ in a general position as
 $\g=\g(n_1,\ldots,n_{\n})=\g(\mb{n})$. In general, for
 $\mb{x},\mb{y}\in\bR^{\n}$ we put $\mb{x}\mb{y}:=\sum^{\n}_{j=1}x_jy_j$ 
 and put $|\mb{x}|:=\sqrt{\mb{x}\mb{x}}=\sqrt{{x_1}^2+\cdots+{x_{\n}}^2}$. 
 Hence $\ell(\g(\mb{n}))=|\mb{n}|$. We may therefore write
 $\c(\g(\mb{n}))=e^{2\pi i\mb{n}\bsym{\a}}$ for some
 $\bsym{\a}=(\a_1,\ldots, \a_{\n})\in \bR^{\n}/\bZ^{\n}$. Thus we
 sometimes identify the unitary character $\c$ of $\G\,(\cong \bZ^{\n})$
 with $\bsym{\a}\in \bR^{\n}/\bZ^{\n}(\cong \hat\bZ^{\n})$ and write
 the Ruelle type $L$-function $L_{\pp}(s,\c)$ as $L_{\pp}(s,\bsym{\a})$:  
\begin{align}
\label{for:def-Lpp}
 L_{\pp}(s,\bsym{\a})=L_{\pp}(s,\bsym{\a};\n)
:=\prod_{\gcd\mb{p}=1}(1-e^{2\pi i\mb{p}\bsym{\a}}e^{-s|\mb{p}|})^{-1}.
\end{align}
 It is easy to see that the Euler product \eqref{for:def-Lpp} converges
 absolutely for $\Re(s)>0$ because
 $\n|\mb{x}|\geq |x_1|+\cdots+|x_{\n}|$ for $\mb{x}\in\bR^{\n}$, whence
 it defines a holomorphic function for $\Re(s)>0$. Also, we notice that
 $L_{\pp}(s,\bsym{\a}(\bsym{\e}))=L_{\pp}(s,\bsym{\a})$ where
 $\bsym{\a}(\bsym{\e})=(\e_1\a_1,\ldots,\e_{\n}\a_{\n})$ for
 $\bsym{\e}=(\e_1,\ldots,\e_{\n}) \in \{\pm\}^{\n}$ and in particular,
 $L_{\pp}(s,\overline\c)=L_{\pp}(s,\c)$. 

 We first show the

\begin{prop}
\label{prop:log-deri-L}
 The logarithmic derivative of $L_{\pp}(s,\bsym{\a};\n)$ is holomorphic
 in $\Re(s)>0$ and has the expression  
\begin{equation}
\label{for:log-deri-L}
 \frac{L'_{\pp}}{L_{\pp}}(s,\bsym{\a};\n)
=2(2\sqrt{\pi})^{\n-1}\G\Bigl(\frac{\n+1}{2}\Bigr)
\sum^{\infty}_{n=1}\frac{\g(n)}{n^{\n+1}}  
\sum_{\mb{m}\in\bZ^{\n}}\frac{\bigl(2\pi|\mb{m}/n+\bsym{\a}|\bigr)^2-\n s^2}
{\Bigl\{s^2+\bigl(2\pi|\mb{m}/n+\bsym{\a}|\bigr)^2\Bigr\}^{\frac{\n+3}{2}}}
\end{equation}
 for $\Re(s)>0$. Here $\g(n)=\prod_{p|n}(1-p)$.
\end{prop}
\begin{proof} 
 The first assertion is clear, because $L_{\pp}(s,\bsym{\a})$ is
 holomorphic and non-zero in $\Re(s)>0$. We introduce now an auxiliary
 function $G(s,\bsym{\a})$ by  
\[
 G(s,\bsym{\a})
:=\prod_{\mb{n}\in \bZ^{\n}\bslO}
(1-e^{2\pi i\mb{n}\bsym{\a}}e^{-s|\mb{n}|})^{-1}.
\]
 Then we have 
\begin{align}
\label{for:logG}
 \log{G(s,\bsym{\a})}
=-\sum_{\mb{n}\in\bZ^{\n}\bslO}\log(1-e^{2\pi i\mb{n}\bsym{\a}}e^{-s|\mb{n}|})
=\sum_{\mb{n}\in \bZ^{\n}\bslO}
\sum^{\infty}_{\ell=1}\frac{1}{\ell}e^{2\pi i\ell\mb{n}\bsym{\a}}
e^{-s\ell|\mb{n}|}.
\end{align}
 It is immediate to see that 
\begin{align*}
\log G(s,\bsym{\a})
=\sum^{\infty}_{d=1}\sum_{\gcd\mb{p}=1}
\sum^{\infty}_{\ell=1}\frac{e^{2\pi i\ell(d\mb{p})\bsym{\a}}
e^{-s\ell|d\mb{p}|}}{\ell}
=\sum^{\infty}_{d=1}\log L_{\pp}(ds,d\bsym{\a}).
\end{align*}
 Hence the M\"obius inversion formula yields
\begin{align}
\label{for:logL}
 \log L_{\pp}(s,\bsym{\a})
=\sum^{\infty}_{m=1}\m(m)\log G(ms,m\bsym{\a}).
\end{align}
 Here $\m(n)$ denotes the M\"obius function defined as
\[
 \m(n)=
\begin{cases}
 1          & \quad \textrm{if $n=1$},\\
 (-1)^k     & \quad \textrm{if $n$ is a product of $k$ distinct primes},\\ 
 0          & \quad \textrm{otherwise}.
\end{cases}
\]
 In fact, since $\sum_{m|n}\m(m)=\d_{n1}$, we see that 
\begin{align*}
  \sum^{\infty}_{m=1}\m(m)\log G(ms,m\bsym{\a})
&=\sum^{\infty}_{m=1}\m(m)\sum^{\infty}_{d=1}\log L_{\pp}(dms,dm\bsym{\a})\\
&=\sum^{\infty}_{n=1}\Big\{\sum_{m|n}\mu(m)\Big\}\log L_{\pp}(ns,n\bsym{\a})
=\log L_{\pp}(s,\bsym{\a}).
\end{align*}

 We now put  
\[
 g(s,\bsym{\a})
:=\sum_{\mb{n}\in\bZ^{\n}\bslO}e^{2\pi i\mb{n}\bsym{\a}}e^{-s|\mb{n}|}. 
\]
 Then by \eqref{for:logG} we easily get the relation
\begin{align}
\label{for:G-and-g}
 \log G(s,\bsym{\a})
=\sum^{\infty}_{\ell=1}\frac{1}{\ell}g(\ell s,\ell\bsym{\a}).
\end{align}
 Let $f(\mb{x})$ be a smooth rapidly decreasing function on $\bR^{\n}$.
 Write the Fourier transform of $f(\mb{x})$ by
 $\hat{f}(\mb{y})=\int_{\bR^{\n}}f(\mb{x})e^{2\pi i\mb{y}\mb{x}}d\mb{x}$.
 It can be calculated (see, e.g., \cite{Muller1998}) as    
\[
 (e^{-2\pi t|\mb{x}|}e^{-2\pi i\mb{z}\mb{x}})^{\wedge}(\mb{y})
=\frac{2t}{\Area(S^{\n})}\cdot
 \bigl(t^2+|\mb{y}-\mb{z}|^2\bigr)^{-\frac{\n+1}{2}}
 \qquad (\Re(t)>0),
\]
 where
 $\Area(S^{\n})=2{\pi^{\frac{\n+1}{2}}}/{\G\bigl(\frac{\n+1}{2}\bigr)}$. 
 Therfore, the Poisson summation formula applying to the function
 $e^{-2\pi t|\mb{x}|}e^{2\pi i\bsym{\a}\mb{x}}$ reads  
\begin{align}
\label{for:PSF}
 \sum_{\mb{n}\in\bZ^{\n}}e^{-2\pi t|\mb{n}|}e^{-2\pi i\mb{z}\mb{n}} 
=\frac{2t}{\Area(S^{\n})}\sum_{\mb{m}\in\bZ^{\n}}
\bigl(t^2+|\mb{m}-\mb{z}|^2\bigr)^{-\frac{\n+1}{2}}.
\end{align}
 This shows that for $\Re(s)>0$ we have 
\[
 g(s,\bsym{\a})+1
=\frac{2(2\pi)^{\n} s}{\Area(S^{\n})}\sum_{\mb{m}\in\bZ^{\n}}
\Big\{s^2+\bigl(2\pi|\mb{m}+\bsym{\a}|\bigr)^2\Big\}^{-\frac{\n+1}{2}}.
\]
 Differentiating this equation with respect to $s$, we obtain
\begin{align}
\label{for:diff-g}
 g'(s,\bsym{\a})
=\frac{2(2\pi)^{\n}}{\Area(S^{\n})}\sum_{\mb{m}\in\bZ^{\n}}\frac{\bigl(2\pi|\mb{m}+\bsym{\a}|\bigr)^2-\n s^2} 
{\Big\{s^2+\bigl(2\pi|\mb{m}+\bsym{\a}|\bigr)^2\Big\}^{\frac{\n+3}{2}}}.
\end{align}
 By the relations \eqref{for:logL} and \eqref{for:G-and-g}, we have 
\[
 \frac{L'_{\pp}}{L_{\pp}}(s,\bsym{\a})
=\sum^{\infty}_{m=1} m\m(m)\frac{G'}{G}(ms,m\bsym{\a}),\quad
\frac{G'}{G}(s,\bsym{\a})
=\sum^{\infty}_{\ell=1} g'(\ell s,\ell\bsym{\a}).
\]
 From these equations it follows that 
\begin{equation}
\label{for:g-log-deri-L}
 \frac{L'_{\pp}}{L_{\pp}}(s,\bsym{\a})
=\sum^{\infty}_{m=1}\sum^{\infty}_{\ell=1} m\m(m)g'(m\ell s,m\ell\bsym{\a})
=\sum^{\infty}_{n=1}\g(n)g'(ns,n\bsym{\a}).
\end{equation}
 Here we have used an easily verified formula $\sum_{m|n}m\m(m)=\g(n)$.
 Hence, it follows from \eqref{for:diff-g} and \eqref{for:g-log-deri-L}
 that    
\begin{align*}
 \frac{L'_{\pp}}{L_{\pp}}(s,\bsym{\a})
&=\frac{2(2\pi)^{\n}}{\Area(S^{\n})}
\sum^{\infty}_{n=1}\g(n)\sum_{\mb{m}\in\bZ^{\n}}
\frac{\bigl(2\pi|\mb{m}+n\bsym{\a}|\bigr)^2-\n (ns)^2}{\Big\{(ns)^2+\bigl(2\pi|\mb{m}+n\bsym{\a}|\bigr)^2\Big\}^{\frac{\n+3}{2}}}\\
&=2(2\sqrt{\pi})^{\n-1}\G\Bigl(\frac{\n+1}{2}\Bigr)
\sum^{\infty}_{n=1}\frac{\g(n)}{n^{\n+1}} 
\sum_{\mb{m}\in\bZ^{\n}}\frac{\bigl(2\pi|\mb{m}/n+\bsym{\a}|\bigr)^2-\n s^2}
{\Bigl\{s^2+\bigl(2\pi|\mb{m}/n+\bsym{\a}|\bigr)^2\Bigr\}^{\frac{\n+3}{2}}}.  
\end{align*}
 This shows the proposition. 
\end{proof}

 Put
\begin{align*}
 C(\n):
&=2(2\sqrt{\pi})^{\n-1}\G\Bigl(\frac{\n+1}{2}\Bigr),\\
 \Phi(s,\bsym{\a},t;\n):
&=\sum_{n\in\bN,\ \mb{m}\in\bZ^{\n}}\frac{\g(n)}{n^{\n+1}}
\Bigl(s^2+\bigl(2\pi|\mb{m}/n+\bsym{\a}|\bigr)^2\Bigr)^{-t}.
\end{align*}
 Then the formula \eqref{for:log-deri-L} can be written as 
\begin{equation}
\label{for:Phi-log-deri-L}
 \frac{L'_{\pp}}{L_{\pp}}(s,\bsym{\a};\n)
=C(\n)\biggl(\Phi\Bigl(s,\bsym{\a},\frac{\n+1}{2};\n\Bigr)
-(\n+1)s^2\Phi\Bigl(s,\bsym{\a},\frac{\n+3}{2};\n\Bigr)\biggr).
\end{equation}
 
 The following theorem shows that $\z_{\pp}(s;\n):=L_{\pp}(s,\mb{0};\n)$
 can not be extended meromorphically across the imaginary axis $\Re(s)=0$ if
 $\n=2,4$ or $8$, that is, the cases where the number $r_{\n}(n)$ of
 expressions of a positive integer $n$ by a sum of $\n$ integers square
 is essentially multiplicative.

\begin{thm}
\label{thm:nb}
 Let $\n=2,4$ or $8$. Then the imaginary axis is a natural boundary of
 $\z_{\pp}(s;\n)$.  
\end{thm}

 From the expression \eqref{for:Phi-log-deri-L}, to prove \Thm{nb}, it
 is sufficient to show that the imaginary axis $\Re(s)=0$ is a natural
 boundary of the function $\Phi(s,\mb{0},t;\n)$. Let
 $\bS_{\n}:=\bigl\{\,\frac{|\mb{m}|}{n}\bigl|\,n\in\bN,\ \mb{m}\in\bZ^{\n}\bigr\}$.
 For a given $x_0\in\bR$, let us denote $R_{\n}(x_0)$ as the coefficient
 of $\bigl(s^2+(2{\pi}x_0)^2\bigr)^{-t}$ in the series
 $\Phi(s,\mb{0},t;\n)$:   
\begin{equation}
\label{for:Phi-coeff}
 \Phi(s,\mb{0},t;\n)
=\sum_{x_0\in\bS_{\n}}R_{\n}(x_0)\bigl(s^2+(2{\pi}x_0)^2\bigr)^{-t}.
\end{equation}
 Note that $R_{\n}(x_0)\equiv 0$ if $x_0\notin \bS_{\n}$.

\begin{lem}
\label{lem:coeff}
 Let $\w{m}$ and $\w{n}$ be positive integers satisfying
 $(\w{m},\w{n})=1$. Then we have
\[
 R_{\n}\bigl(\frac{\w{m}}{\w{n}}\bigr)
=\frac{1}{{\w{n}}^{\n+1}}\sum^{\infty}_{k=1}
\frac{\g(k\w{n})r_{\n}(k^2{\w{m}}^2)}{k^{\n+1}},
\]
 where
 $r_{\n}(n):=\#\bigl\{(m_1,\ldots,m_{\n})\in\bZ^{\n}\,\bigl|\,{m_1}^2+\cdots+{m_{\n}}^2=n\bigr\}$.
\end{lem}
\begin{proof}
 Since $(\w{m},\w{n})=1$, the term
 $\bigl(s^2+(2{\pi}\frac{\w{m}}{\w{n}})^2\bigr)^{-t}$ appears in
 $\Phi(s,\mb{0},t;\n)$ when $n$ and $\mb{m}$ can be written as
 $n=k\w{n}$ and $|\mb{m}|=k\w{m}$ with some $k\in\bN$. Therefore, one
 can calculate as 
\begin{align*}
  R_{\n}\bigl(\frac{\w{m}}{\w{n}}\bigr)
&=\sum^{\infty}_{k=1}\sum_{|\mb{m}|=k\w{m}}\frac{\g(k\w{n})}{(k\w{n})^{\n+1}}
=\frac{1}{{\w{n}}^{\n+1}}\sum^{\infty}_{k=1}
\sum_{m_1,\ldots,m_{\n}\in\bZ \atop {m_1}^2+\cdots+{m_{\n}}^2=k^2{\w{m}}^2}
\frac{\g(k\w{n})}{k^{\n+1}}\\
&=\frac{1}{{\w{n}}^{\n+1}}\sum^{\infty}_{k=1}
\frac{\g(k\w{n})r_{\n}(k^2{\w{m}}^2)}{k^{\n+1}}.
\end{align*}
 This shows the assertion.
\end{proof}

 The following lemma is crucial for the proof of \Thm{nb}.

\begin{lem}
\label{lem:key}
 Assume $\n=2,4$ or $8$. For any prime $p$ and a non-negative integer
 $e$, it holds that 
\begin{equation}
\label{for:key}
 \sum^{\infty}_{n=1}r_{\n}(p^{2(n+e)})p^{-n(\n+1)}
\ne\frac{r_{\n}(p^{2e})}{p-1}.
\end{equation}
\end{lem}
\begin{proof}
 It is a classical result that the arithmetic function
 $\frac{1}{2\n}r_{\n}(n)$ is multiplicative, that is,
 $\frac{1}{2\n}r_{\n}(mn)=\frac{1}{2\n}r_{\n}(m)\frac{1}{2\n}r_{\n}(n)$
 for $(m,n)=1$ if and only if $\n=2,4$ or $8$, and is explicitely given
 (by Jacobi when $\n=4$) as  
\begin{equation}
\label{for:rn}
 \frac{1}{2\n}r_{\n}(n)=
\begin{cases}
 \displaystyle{\sum_{m|n}\c_{-4}(m)} & \quad \textrm{if $\n=2$},\\[12pt]
 \displaystyle{\sum_{m|n,4\nmid m}m} & \quad \textrm{if $\n=4$},\\[12pt]
 \displaystyle{(-1)^n\sum_{m|n}(-1)^m m^3} & \quad \textrm{if $\n=8$},
\end{cases}
\end{equation} 
 where $\c_{-4}(n)$ is the primitive Dirichlet character modulo $4$
 (see, e.g., \cite{BorweinChoi2003}). We check the formula
 \eqref{for:key} for each $\n$.
\begin{enumerate}
 \item The case $\n=2$: Let $\ell$ be a positive integer and $p$ a
       prime. From \eqref{for:rn}, it holds that 
\begin{equation}
\label{for:r2}
 r_{2}(1)=4, \quad 
 r_{2}(2^{\ell})=4  \quad \textrm{and} \quad 
 r_{2}(p^{\ell})=
\begin{cases}
 4(\ell+1) & \textrm{if $p\equiv 1 \pmod 4$,}\\[5pt]
 \begin{cases}
   4 & \quad \textrm{if $\ell$ is even},\\
   0 & \quad \textrm{if $\ell$ is odd},   
 \end{cases}
           & \textrm{if $p\equiv 3 \pmod 4$.}
\end{cases}
\end{equation}
       Therfore we have 
\begin{enumerate}
 \item for $p=2$
\begin{align*}
 \textrm{L.H.S of \eqref{for:key}}
&=4\sum^{\infty}_{n=1}2^{-3n}=\frac{4}{7},\qquad\qquad\qquad\qquad\qquad\qquad\qquad\quad\ \ \\
 \textrm{R.H.S of \eqref{for:key}}
&=4.
\end{align*}
       Hence the claim \eqref{for:key} is true.

 \item for $p\equiv 1 \pmod 4$
\begin{align*}
 \textrm{L.H.S of \eqref{for:key}}
&=\sum^{\infty}_{n=1}4\bigl(2(n+e)+1\bigr)2^{-3n}
 =\frac{8p^3+4(2e+1)(p^3-1)}{(p^3-1)^2},\\
 \textrm{R.H.S of \eqref{for:key}}
&=\frac{4(2e+1)}{p-1}.
\end{align*}
       In this case, since $p^{\ell}\equiv 1\pmod 4$ for all $\ell\ge 1$,
       we obtain \eqref{for:key}. 

 \item for $p\equiv 3 \pmod 4$
\begin{align*}
 \textrm{L.H.S of \eqref{for:key}}
&=4\sum^{\infty}_{n=1}p^{-3n}=\frac{4}{p^3-1},\qquad\qquad\qquad\qquad\qquad\qquad\quad\
 \\
 \textrm{R.H.S of \eqref{for:key}}
&=\frac{4}{p-1}.
\end{align*}
       Hence \eqref{for:key} follows clearly.
\end{enumerate}
 \item The case $\n=4$: It follows from \eqref{for:rn} again that
\begin{equation}
\label{for:r4}
 r_{4}(1)=8, \quad
 r_{4}(2^{\ell})=24
\quad \textrm{and} \quad
 r_{4}(p^{\ell})=8\cdot\frac{p^{\ell+1}-1}{p-1} \quad (p\ne 2).
\end{equation}
       Hence we have
\begin{enumerate}
 \item for $p=2$
\begin{align*}
 \textrm{L.H.S of \eqref{for:key}}
&=24\sum^{\infty}_{n=1}2^{-5n}=\frac{24}{31},\qquad\qquad\qquad\qquad\qquad\qquad\qquad\
 \ \\
 \textrm{R.H.S of \eqref{for:key}}
&=\begin{cases}
   8 \quad & \textrm{if} \quad e=0,\\
   24 \quad & \textrm{if} \quad e>0.  
  \end{cases}
\end{align*}
       This shows that \eqref{for:key} is true.

 \item for $p\ne 2$
\begin{align*}
 \textrm{L.H.S of \eqref{for:key}}
&=8\sum^{\infty}_{n=1}\frac{p^{2(n+e)+1}-1}{p-1}p^{-5n}
=\frac{8}{p-1}\Bigl(\frac{p^{2e+1}}{p^3-1}-\frac{1}{p^5-1}\Bigr),\\
 \textrm{R.H.S of \eqref{for:key}}
&=\frac{8}{p-1}\cdot \frac{p^{2e+1}-1}{p-1}.
\end{align*}
       Since $p>3$, we have
\[
 \frac{p^{2e+1}-1}{p-1}-\Bigl(\frac{p^{2e+1}}{p^3-1}-\frac{1}{p^5-1}\Bigr)
>\frac{p^{2e+1}-1}{p-1}-\frac{p^{2e+1}}{p^3-1}
>0.
\]
       Hence the claim \eqref{for:key} follows.
\end{enumerate}
 \item The case $\n=8$: We have
\begin{equation}
\label{for:r8}
 r_{8}(1)=16, \quad
 r_{8}(2^{\ell})=16\cdot\frac{2^{3(\ell+1)}-15}{7}
\quad \textrm{and} \quad
 r_{8}(p^{\ell})=16\cdot\frac{p^{3(\ell+1)}-1}{p^3-1} \quad (p\ne 2).
\end{equation}
       Hence it holds that 
\begin{enumerate}
 \item for $p=2$
\begin{align*}
 \textrm{L.H.S of \eqref{for:key}}
&=16\sum^{\infty}_{n=1}\frac{2^{3(2(n+e)+1)}-15}{7}2^{-9n}
=\frac{16}{7}\Bigl(\frac{2^{3(2e+1)}}{7}-\frac{15}{511}\Bigr),\qquad\\
 \textrm{R.H.S of \eqref{for:key}}
&=\begin{cases}
   16 & \quad \textrm{if $e=0$,}\\[5pt]
   \displaystyle{16\cdot\frac{2^{3(2e+1)}-15}{7}} & \quad \textrm{if $e>0$.} 
  \end{cases}
\end{align*}
       Hence the claim \eqref{for:key} follows.
 
 \item for $p\ne 2$
\begin{align*}
 \textrm{L.H.S of \eqref{for:key}}
&=16\sum^{\infty}_{n=1}\frac{p^{3(2(n+e)+1)}-1}{p^3-1}p^{-9n}
=\frac{16}{p^3-1}\Bigl(\frac{p^{3(2e+1)}}{p^3-1}-\frac{1}{p^9-1}\Bigr),\\
 \textrm{R.H.S of \eqref{for:key}}
&=\frac{16}{p-1}\cdot\frac{p^{3(2e+1)}-1}{p^3-1}.
\end{align*}
       By the same argument in the cases of $\n=4$ and $p\ne 2$, one can
       show \eqref{for:key}.      
\end{enumerate}
\end{enumerate}
 This completes the proof of the lemma.
\end{proof}

 We now prove \Thm{nb}.

\begin{proof}
[Proof of \Thm{nb}]
 Since $\bQ$ is dense in $\bR$, it is sufficient to show that
 $R_{\n}(\xi)\ne 0$ for any positive $\xi\in\bQ$ by the expression
 \eqref{for:Phi-coeff} of $\Phi(s,\bsym{\a},t;\n)$. Let $\w{m}$ and
 $\w{n}$ be positive integers satisfying $(\w{m},\w{n})=1$. Let
 $\w{m}={p_1}^{e_1}\cdots{p_s}^{e_s}$ and
 $\w{n}={q_1}^{f_1}\cdots{q_t}^{f_t}$ be their prime
 factorizations. Note that $(p_i,q_j)=1$ for any $i$ and $j$. We shall 
 prove $R_{\n}\bigl(\frac{\w{m}}{\w{n}}\bigr)\ne 0$. Since the functions
 $\g(n)$ and $\frac{1}{2\n}r_{\n}(n)$ with $\n=2,4$ or $8$ are
 multiplicative, we have from \Lem{coeff} that  
\begin{align*}
 {\w{n}}^{\n+1}R_{\n}\bigl(\frac{\w{m}}{\w{n}}\bigr)
&=\sum_{m_1,\ldots,m_t\ge 0 \atop n_1,\ldots,n_s\ge 0}
\sum_{(k,p_i)=1 \atop (k,q_j)=1}
\frac{\g\bigl(k\prod^{t}_{j=1}{q_j}^{m_j+f_j}\prod^{s}_{i=1}{p_i}^{n_i}\bigr)r_n\bigl(k^2\prod^{t}_{j=1}{q_j}^{2m_j}\prod^{s}_{i=1}{p_i}^{2(n_i+e_i)}\bigr)}{\bigl(k\prod^{t}_{j=1}{q_j}^{m_j}\prod^{s}_{i=1}{p_i}^{n_i}\bigr)^{\n+1}}\\
&=\biggl(\prod^{t}_{j=1}\frac{1-q_j}{2\n}\sum^{\infty}_{m_j=0}\frac{r_{\n}({q_j}^{2m_j})}{{q_j}^{{m_j}(\n+1)}}\biggr)
\biggl(\prod^{s}_{i=1}\frac{1}{2\n}\sum^{\infty}_{n_i=0}\frac{\g({p_i}^{n_i})r_{\n}({p_i}^{2(n_i+e_i)})}{{p_i}^{{n_i}(\n+1)}}\biggr)
\biggl(\sum_{(k,p_i)=1 \atop
 (k,q_j)=1}\frac{\g(k)r_{\n}(k^2)}{k^{\n+1}}\biggr)\\ 
&=\prod^{t}_{j=1}\frac{1-q_j}{2\n}E_{\n,q_j}\bigl(\frac{\w{m}}{\w{n}}\bigr)
\times\prod^{s}_{i=1}\frac{1}{2\n}F_{\n,p_i}\bigl(\frac{\w{m}}{\w{n}}\bigr)
\times G_{\n}\bigl(\frac{\w{m}}{\w{n}}\bigr),
\end{align*}
 where 
\begin{align*}
 E_{\n,q_j}\bigl(\frac{\w{m}}{\w{n}}\bigr):
&=\sum^{\infty}_{m_j=0}\frac{r_{\n}({q_j}^{2m_j})}{{q_j}^{{m_j}(\n+1)}}, \quad
 F_{\n,p_i}\bigl(\frac{\w{m}}{\w{n}}\bigr):
=\sum^{\infty}_{n_i=0}\frac{\g({p_i}^{\n_i})r_{\n}({p_i}^{2(n_i+e_i)})}{{p_i}^{{n_i}(\n+1)}}
\quad \textrm{and} \quad\\
 G_{\n}\bigl(\frac{\w{m}}{\w{n}}\bigr):
&=\sum_{(k,p_i)=1 \atop (k,q_j)=1}\frac{\g(k)r_{\n}(k^2)}{k^{\n+1}}.
\end{align*}
 Note that $\g({q_j}^{m_j+f_j})=\g(q_j)=1-q_j$, because $f_j\ge 1$. Since
 $E_{\n,q_j}\bigl(\frac{\w{m}}{\w{n}}\bigr)\ne 0$ for all $j$, it
 suffices to show that $F_{\n,p_i}(\frac{\w{m}}{\w{n}})\ne 0$ for all
 $i$ and $G_{\n}\bigl(\frac{\w{m}}{\w{n}}\bigr)\ne 0$. 
\begin{enumerate}
 \item The function $F_{\n,p_i}\bigl(\frac{\w{m}}{\w{n}}\bigr)$: Since
\begin{align*}
 F_{\n,p_i}\bigl(\frac{\w{m}}{\w{n}}\bigr)
&=r_{\n}(p^{2e_i})-(p_i-1)\sum^{\infty}_{n_i=1}
\frac{r_{\n}({p_i}^{2(n_i+e_i)})}{{p_i}^{{n_i}(\n+1)}},
\end{align*}
       we have
\[
  F_{\n,p_i}\bigl(\frac{\w{m}}{\w{n}}\bigr)\ne 0 \iff
 \sum^{\infty}_{n_i=1}r_{\n}({p_i}^{2(n_i+e_i)}){p_i}^{-{n_i}(\n+1)}
 \ne \frac{r_{\n}(p^{2e_i})}{p_i-1}.
\]
       By \Lem{key} (the case of $e=e_i\ge 1$), the right-above is true,
       whence the claim follows.  
 
 \item The function $G_{\n}(\frac{\w{m}}{\w{n}})$: Since $\g(n)$ and 
       $\frac{1}{2\n}r_{\n}(n^2)$ are multiplicative, the series
\[
 G_{\n}:=\sum^{\infty}_{k=1}\frac{\g(k)r_{\n}(k^2)}{k^{\n+1}}
\]
       has an Euler product expression. In fact, we have
       $G_{\n}=2\n\prod_{p:\textrm{prime}}G_{\n,p}$, where
\[
  G_{\n,p}:
=\frac{1}{2\n}\sum^{\infty}_{\ell=0}\g(p^{\ell})r_{\n}(p^{2\ell})
p^{-{\ell}(\n+1)}
=1-\frac{p-1}{2\n}\sum^{\infty}_{\ell=1}r_{\n}(p^{2\ell})p^{-{\ell}(\n+1)}
\ne 0.
\]
       The fact $G_{\n,p}\ne 0$ is, in fact, derived from \Lem{key}
       (i.e., look at the case $e=0$). It follows that  
\[
 G_{\n}(\frac{\w{m}}{\w{n}})=\prod_{p:\textrm{prime} \atop p\ne p_i, q_j} G_{\n,p}\ne 0.
\]
\end{enumerate}
 This completes the proof of the theorem.
\end{proof}

 Since, for a given $x_0\in\bR$ and any $\e>0$, there exists
 $\mb{m}\in\bZ^{\n}$ and $n\in\bN$ satisfying
 $\Big|2\pi\big|\frac{\mb{m}}{n}+\bsym{\a}\big|-x_0\Big|<\e$, we
 may expect that the set of singularities of
 $\Phi(s,\bsym{\a},t;\n)$ forms a dense subset of $\bR$. In other words,
 this expectation suggests the

\begin{conj}
 The imaginary axis is a natural boundary of $L_{\pp}(s,\bsym{\a};\n)$. 
\end{conj}

\section{Determinants of Laplacians}

 As we have shown, the Ruelle type $L$-function  $L_{\pp}(s,\c)$
 (defined by the Euler product) does not seem to be extended to the
 entire plane $\bC$ as a meromorphic function. Hence, in particular, it
 may not have a determinant expression for a ``Laplacian''. Thus,
 contrary to the discussion above, one may naturally ask the question;
 does the function defined by a determinant of the Laplacian on the
 torus $\bT^{\n}:=\bR^{\n}/\bZ^{\n}$ have any Euler product expression
 in an appropriate sense ?    

\smallskip

 We now recall the definition of the determinant of a Laplacian via zeta
 regularization. When the sequence $\{\l_n\}_{n=1,2,\ldots}$ consisting
 of the eigenvalues of some (infinite dimensional) operator $\Om$ is
 given, we can define the determinant of $\Om$ by
 $\det(\Om):=\regprod^{\infty}_{n=1}\l_n$. Here, in general, the zeta
 regularized product $\regprod^{\infty}_{n=1} a_n$ of the sequence
 $\{a_n\}_{n=1,2,\ldots}$ is defined by 
 $\regprod^{\infty}_{n=1}a_n:=\exp\bigl(-\frac{d}{ds}\sum^{\infty}_{n=1}a_n^{-s}\bigl|_{s=0}\bigr)$
 when the Dirichlet series $\sum^{\infty}_{n=1} a_n^{-s}$ converges
 absolutely for $\Re(s) \gg 1$ and can be analytically continued to a
 holomorphic function around the origin $s=0$ (see \cite{Deninger1991}
 and also \cite{ZR}).  

 Concerning the question above, our observation is as follows: Let
 $\bsym{\a}\in\hat\bZ^{\n}\cong\bR^{\n}/\bZ^{\n}$. Let
 $\D_{\n,\bsym{\a}}:=-\frac1{4{\pi}^2}\sum^{\n}_{j=1}\frac{\p^2}{\p x_j^2}$
 be the Laplacian acting on the space of smooth sections of the line
 bundle defined by the unitary character
 $\c_{\bsym{\a}}:\bZ^{\n}\ni\mb{n} \to e^{2\pi i\mb{n}\bsym{\a}}\in U(1)$
 over the torus $\bT^{\n}$. It is easy to see that the $L^2$-eigenvalue of
 $\D_{\n,\bsym{\a}}$ is of the form $|\mb{m}+\bsym{\a}|^2$. Hence we have
\[
 \det(\D_{\n,\bsym{\a}}+s^2)
=\regprod_{\mb{m}\in \bZ^{\n}}\bigl(|\mb{m}+\bsym{\a}|^2+s^2\bigr).
\]  
 We now show the following theorem. To simplify the description, we
 write $Z_1(s) \cong_k Z_2(s)$ for two functions $Z_i(s)$ when there is
 a polynomial $P(s)$ of degree $k$ satisfying $Z_1(s)=Z_2(s)e^{P(s)}$.

\begin{thm}
\label{thm:DetLaplacian}
 Let $\bsym{\a}=(\a_1,\ldots,\a_{\n})\in\hat\bZ^{\n}=\bR^{\n}/\bZ^{\n}$.
 Assume $\Re(s)>0$. 
\begin{enumerate}
 \item Suppose that $\n=2\ell+1$ $(\ell \ge 0)$. Then
\begin{align}
\label{for:det-odd}
&\det(\D_{2\ell+1,\bsym{\a}}+s^2)
=\regprod_{\mb{m}\in \bZ^{2\ell+1}}\bigl(|\mb{m}+\bsym{\a}|^2+s^2\bigr)\\
&\qquad \cong_{2\ell} 
\exp\biggl(-\frac{(-2\pi)^{\ell+1}}{(2\ell+1)!!} s^{2\ell+1}-
\sum_{\mb{n} \in \bZ^{2\ell+1}\bslo}\frac{1}{|\mb{n}|^{\ell+1}}
\sum^{\ell}_{k=0}\frac{c^{(\ell)}_k}{(2\pi|\mb{n}|)^k}s^{\ell-k}
e^{2\pi i\mb{n}\bsym{\a}}e^{-2\pi|\mb{n}|s}\biggr),\nonumber
\end{align}
 where the sequence $c^{(\ell)}_k$ $(0\le k \le \ell)$ is defined by
\begin{equation}
\label{for:cl}
  c^{(\ell)}_0=1, \qquad c^{(\ell)}_k
=\begin{cases}
 1 & \quad \textrm{if $\ell=0$},\\
 \displaystyle{\frac{1}{{2^k}{k!}}\prod^k_{j=1-k}(\ell+j)} & \quad \textrm{if $l\ge 1$},
 \end{cases}
\quad (k\ge 1).
\end{equation}
 Moreover, when $\n=1$ $(\textrm{that is},\ \ell=0)$, we have 
\begin{equation}
\label{for:det-odd-1}
 \det(\D_{1,\a}+s^2)
=\regprod^{\infty}_{m=-\infty}\big((m+\a)^2+s^2\big)
=\exp\Bigl(2\pi s-2\sum^{\infty}_{n=1}
 \frac{\cos(2\pi n\a)}{n}e^{-2\pi ns}\Bigr).
\end{equation}
 \item Suppose that $\n=2\ell$ $(\ell \ge 1)$. Then
\begin{align}
\label{for:det-even}
&\det(\D_{2\ell,\bsym{\a}}+s^2)
=\regprod_{\mb{m}\in \bZ^{2\ell}}\bigl(|\mb{m}+\bsym{\a}|^2+s^2\bigr)\\
&\qquad \cong_{2\ell} 
\exp\biggl(\frac{4(-1)^{\ell}{\pi}^{\ell}}{(\ell-1)!}s^{2\ell}\log{s}-4{\ell}s^{\ell}
\sum_{\mb{n}\in \bZ^{2\ell}\bslo}\frac{1}{|\mb{n}|^{\ell}}e^{2\pi i\mb{n}\bsym{\a}}K_{\ell}(2\pi|\mb{n}|s)\biggr),\qquad\qquad\quad\nonumber
\end{align}
 where $K_{\ell}(s)$ is the $K$-Bessel function.
\end{enumerate}
\end{thm}

 Before starting the proof, we remark that the formula in the theorem
 can be considered as a sort of an Euler product expression of the
 determinant of the Laplacian. Actually, when $\n=1$, from the identity
 \eqref{for:det-odd-1} we easily find that
\begin{align}
\label{for:LerchGen}
 \det(\D_{1,\a}+s^2)
=e^{2\pi s}(1-e^{2\pi i \a}e^{-2\pi s})(1-e^{-2\pi i\a}
 e^{-2\pi s}).  
\end{align}
 Since the only primitive closed geodesic on $S^1=\bR/\bZ$ is the circle
 of length $2\pi$, the right hand side can be regarded as the Euler
 product (over the primitive element of $\bZ$). Hence it is reasonable
 to regard the right hand side of
 \eqref{for:det-odd}/\eqref{for:det-even} as an Euler product.   

\smallbreak

 We first recall the Bessel functions. Let $J_{\l}(x)$ be the Bessel
 function of the first kind of order $\l$. Namely, $J_{\l}(x)$ satisfies
 the Bessel equation of order $\l$:  
\begin{equation}
\label{Bessel-equation}
 \biggl(\frac{d^2}{dx^2}+\frac{1}{x}\frac{d}{dx}
+\Bigl(1-\frac{{\l}^2}{x^2}\Bigr)\biggr)J_{\l}(x)=0. 
\end{equation}
 Then the $K$-Bessel function $K_{\l}(x)$ is defined as  
\begin{equation}
\label{def:K-Bessel}
 K_{\l}(x):=\frac{\pi}{2}\frac{I_{-\l}(x)-I_{\l}(x)}{\sin{\pi\l}},
 \qquad I_{\l}(x):=e^{-\l\pi i/2}J_{\l}(ix).
\end{equation}
 The following properties are well known (See, e.g.,
 \cite{AndrewsAskeyRoy1999}).

\begin{lem}
\label{lem:Bessel}
 $(i)$\ We have
\begin{equation}
\label{for:J0}
 \int^{2\pi}_{0}e^{ix\cos{\th}}d{\theta}=2\pi J_0(x).
\end{equation}

 $(ii)$\ For $\Re(\xi)>-1$ and $\Re(\eta)>-1$, we have
\begin{equation}
\label{for:Jmn}
 \frac{x^{\xi}y^{\eta}J_{\x+\eta+1}((x^2+y^2)^{1/2})}{(x^2+y^2)^{(\xi+\eta+1)/2}}
=\int^{\pi/2}_{0}J_{\xi}(x\sin{\th})J_{\eta}(y\cos{\th})
\sin^{\xi+1}{\th}\cos^{\eta+1}{\th}d{\th}.
\end{equation}

 $(iii)$\ For $a,b>0$ and $-1<\Re(\l)<2\Re(\m)+3/2$, we have
\begin{equation}
\label{for:Bessel-int}
\int^{\infty}_{0}\frac{x^{\l+1}J_{\l}(bx)}{(x^2+a^2)^{\m+1}}dx
=\frac{a^{\l-\m}b^{\m}}{2^{\l}\G(\l+1)}K_{\l-\m}(ab), 
\end{equation}\qed
\end{lem}

\smallbreak
\begin{proof}
[Proof of Theorem\,\ref{thm:DetLaplacian}]
 By the definition of the regularized product, for a sufficiently large
 integer $j$, we have
\begin{align*}
&\ \ \ \frac{\p^j}{\p t^j}
\log\Bigl(\regprod_{\mb{m}\in \bZ^{\n}}\bigl(|\mb{m}+\bsym{\a}|^2+t\bigr)\Bigr)\\
&=-\frac{d}{dw}
\Big[(-1)^j w(w+1)\cdots(w+j-1)
\sum_{\mb{m}\in \bZ^{\n}}\bigl(|\mb{m}+\bsym{\a}|^2+t\bigr)^{-w-j}
\Big]\Bigg|_{w=0}\\
&=(-1)^{j-1}(j-1)!
\sum_{\mb{m}\in \bZ^{\n}}\bigl(|\mb{m}+\bsym{\a}|^2+t\bigr)^{-j}.
\end{align*}
 Put $j=\ell+1$ for both cases $\n=2\ell+1$ and $2\ell$. Since
 $j=\ell+1>\n/2$, the series in the right hand side converges
 absolutely. Since $\frac{\p}{\p t}=\frac{1}{2s}\frac{\p}{\p s}$ if
 $t=s^2$, we have     
\begin{align}
\label{for:even-odd}
\Bigl(\frac{1}{2s}\frac{\p}{\p s}\Bigr)^{\ell+1}\log 
\Bigl(\regprod_{\mb{m}\in \bZ^{\n}}
 \bigl(|\mb{m}+\bsym{\a}|^2+s^2\bigr)\Bigr) 
=(-1)^{\ell} {\ell}! 
\sum_{\mb{m}\in \bZ^{\n}}\bigl(|\mb{m}+\bsym{\a}|^2+s^2\bigl)^{-\ell-1}.
\end{align}

\begin{enumerate}
 \item The case $\n=2\ell+1$ $(\ell\ge 0)$: By the Poisson summation
       formula \eqref{for:PSF} with $-\ell-1=-\frac{\n+1}{2}$, we see
       that    
\begin{align}
\label{for:PSF-odd}
\frac{\p}{\p s}
\Bigl(\frac{1}{2s}\frac{\p}{\p s}\Bigr)^{\ell}
\log \Bigl(\regprod_{\mb{m}\in \bZ^{2\ell+1}}
 \bigl(|\mb{m}+\bsym{\a}|^2+s^2\bigr)\Bigr) 
=2(-1)^{\ell} {\pi}^{\ell+1}\sum_{\mb{n}\in \bZ^{2\ell+1}}
e^{2\pi i \mb{n}\bsym{\a}}e^{-2\pi s|\mb{n}|}.
\end{align}
        To integrate this equation, we need the following lemma.

\begin{lem}
 Define a polynomial $P_{\ell}(s,a)\; (\ell=0,1,2,\ldots)$ in $s$ by
\[
 P_{\ell}(s,a):=
\begin{cases}
 \displaystyle{\frac{2^{\ell+1}}{(2\ell+1)!!}s^{2\ell+1}}
  & \textrm{if}\quad a= 0,\\[10pt]
 \displaystyle{\Bigl(-\frac{2}{a}\Bigr)^{\ell+1}\sum^{\ell}_{k=0}
 c^{(\ell)}_ka^{-k}s^{\ell-k}} & \textrm{if}\quad a\ne 0,
\end{cases}
\]   
 where $c^{(\ell)}_k$ are given in \eqref{for:cl}. Then, $P_{\ell}(s,a)$
 satisfies the equation  
\begin{align}
\label{for:ladder}
 \frac{\p}{\p s}\Bigl(\frac{1}{2s}\frac{\p}{\p s}\Bigr)^{\ell}
\big\{e^{-as}P_{\ell}(s,a)\big\}=2e^{-as}.
\end{align}
\end{lem}
\begin{proof}
 To prove \eqref{for:ladder}, it is enough to show
\[
 -aP_{\ell}(s,a)+P_{\ell}'(s,a)
=2sP_{\ell-1}(s,a) \quad (\ell=1,2,\ldots).
\]   
 Hence the case $a=0$ is clear. If $a\ne 0$, this is equivalent to the
 recursion 
\[
  c^{(\ell)}_k-c^{(\ell-1)}_k
=(\ell-k+1)c^{(\ell)}_{k-1} \quad (k=1,2,\ldots,\ell).
\]
 One can easily check this recursion from the definition of
 $c^{(\ell)}_k$, whence the claim follows.
\end{proof}

       From the equation \eqref{for:PSF-odd}, it follows immediately
       that 
\begin{align*} 
\frac{\p}{\p s}
\Bigl(\frac{1}{2s}\frac{\p}{\p s}\Bigr)^{\ell}
& \Bigg[\log 
\Bigl(\regprod_{\mb{m}\in \bZ^{2\ell+1}}
 \bigl(|\mb{m}+\bsym{\a}|^2+s^2\bigr)\Bigr)\\ 
&-(-1)^{\ell} \pi^{\ell+1}
\biggl\{P_{\ell}(s,0)
+\sum_{\mb{n}\in \bZ^{2\ell+1}\bslO}P_{\ell}(s,2\pi|\mb{n}|)
e^{2\pi i \mb{n}\bsym{\a}}e^{-2\pi s|\mb{n}|}\biggr\}\Bigg]=0.
\end{align*}
       Hence the equation \eqref{for:det-odd} follows. In particular,
       when $\n=1$, the equation \eqref{for:det-odd} becomes
\begin{equation*}
 \regprod^{\infty}_{m=-\infty}\bigl((m+\a)^2+s^2\bigr)\cong_0 
\exp\bigl(2\pi s-2\sum^{\infty}_{n=1} \frac{\cos(2\pi n\a)}{n}
e^{-2\pi ns}\bigr). 
\end{equation*}
       In other words, there exists a constant $C(\a)$ such that 
\begin{equation}
\label{for:ca}
 \regprod^{\infty}_{m=-\infty} \bigl((m+\a)^2+s^2\bigr)
=C(\a) e^{2\pi s}(1-e^{2\pi i \a}e^{-2\pi s})(1-e^{-2\pi i\a} e^{-2\pi s}).
\end{equation}
       Since one knows in \cite{GLerch} that
\[
 \regprod^{\infty}_{m=0}\bigr((m+\a)^2+s^2\bigr)
=\frac{2\pi}{\G(\a+is)\G(\a-is)},
\]
       the left hand side of \eqref{for:ca} is calculated as 
\begin{align*}
&\regprod^{\infty}_{m=0}\bigr((m+\a)^2+s^2\bigr)\regprod^{\infty}_{m=0}\bigr((m+1-\a)^2+s^2\bigr)\\
&\qquad\ =4{\pi}^2\bigl\{\G(\a+is)\G(1-(\a+is))\G(\a-is)\G(1-(\a-is))\bigr\}^{-1}\\
&\qquad\ =4\sin{\pi{(\a+is)}}\sin{\pi{(\a-is)}}\\
&\qquad\ =e^{2\pi s}(1-e^{2\pi i \a}e^{-2\pi s})(1-e^{-2\pi i\a} e^{-2\pi s}).
\end{align*}
       Hence it proves $C(\a)=1$.

 \item The case $\n=2\ell$ $(\ell\ge 1)$: We calculate the right hand
       side of \eqref{for:even-odd} for $\n=2\ell$. Notice that 
\begin{equation}
\label{for:Fourier-gl}
 \Bigl\{\bigl(|\mb{x}+\bsym{\a}|^2+s^2\bigl)^{-\ell-1}\Bigr\}^{\wedge}(\mb{y})
=e^{-2\pi i\bsym{\a}\mb{y}}g_{\ell}(\mb{y},s),
\end{equation}
 where
\begin{align*}
 g_{\ell}(\mb{y},s):
&=\int_{\bR^{2\ell}}
\frac{e^{2\pi i\mb{x}\mb{y}}}{\bigl(|\mb{x}|^2+s^2\bigr)^{l+1}}d\mb{x}\\
&=\int_{\bR^{2\ell}}\frac{\exp\bigl(2\pi i\sum^{\ell}_{j=1}(x_{2j-1}y_{2j-1}+x_{2j}y_{2j})\bigr)}{\bigl(\sum^{\ell}_{j=1}({x_{2j-1}}^2+{x_{2j}}^2)+s^2\bigr)^{\ell+1}}\prod^{\ell}_{j=1}dx_{2j-1}dx_{2j}.
\end{align*}

  Introduce polar coordinates in both $\mb{x}, \mb{y}\in\bR^{2\ell}$ by
\[
\begin{cases}
 x_{2j-1}=r_j\cos{{\th{_j}}}, & x_{2j}=r_j\sin{{\th}_j},\\
 y_{2j-1}=R_j\cos{{\psi}_j},  & y_{2j}=R_j\sin{{\psi}_j}
\end{cases}
 \qquad (j=1,\ldots,\ell).
\]
 Then we have  
\begin{align*}
  g_{\ell}((\mb{R},\bsym{\psi}),s)
&=g_{\ell}((R_1\cos{{\psi}_1},R_1\sin{{\psi}_1},\ldots,R_{\ell}\cos{{\psi}_{\ell}},R_{\ell}\sin{{\psi}_{\ell}}),s)\\ 
&=\underbrace{\int^{\infty}_{0}\cdots\int^{\infty}_{0}}_{\ell}\frac{\prod^{\ell}_{j=1}r_jdr_j}{\bigl(\sum^{\ell}_{j=1}{r_j}^2+s^2\bigr)^{\ell+1}}\prod^{\ell}_{j=1}\biggl(\int^{2\pi}_{0}\exp\bigl(2\pi
 ir_jR_j\cos{({\th}_j-{\psi}_j)}\bigr)d{\th}_j\biggr)\\
&=(2\pi)^{\ell}\underbrace{\int^{\infty}_{0}\cdots\int^{\infty}_{0}}_{\ell}\frac{\prod^{\ell}_{j=1}r_jJ_0(2\pi r_jR_j)dr_j}{\bigl(\sum^{\ell}_{j=1}{r_j}^2+s^2\bigr)^{\ell+1}}.
\end{align*}
       In the last equality, we have used the formula
       \eqref{for:J0}. Now it is clear that the integral
       $g_{\ell}((\mb{R},\bsym{\psi}),s)$ does not depend on the
       variable $\bsym{\psi}=(\psi_1,\ldots,\psi_{\ell})$. Hence we
       write $g_{\ell}((\mb{R},\bsym{\psi}),s)$ as $g_{\ell}(\mb{R},s)$.

\begin{lem}
\label{lem:gl}
 For $\ell \ge 1$, we have 
\begin{equation}
\label{for:glR}
 g_{\ell}(\mb{R},s)
=\frac{4{\pi}^{\ell+1}|\mb{R}|s^{-1}}{(\ell-1)!}K_1(2\pi|\mb{R}|s),
\end{equation}
 where $|\mb{R}|=\sqrt{{R_1}^2+\cdots+{R_{\ell}}^2}$.
\end{lem}
\begin{proof}
 Put $r_j=r\bigl(\prod^{j-1}_{k=1}\sin{{\th}_k}\bigr)\cos{{\th}_{j}}$
 for $j=1,\ldots,\ell$ and ${\th}_{\ell}=0$. Then one knows
\begin{align*} 
 \sum^{\ell}_{j=1}{r_j}^2=r^2, \qquad
 \prod^{\ell}_{j=1}r_j=r^{\ell}\prod^{\ell-1}_{j=1}\sin^{\ell-j}{{\th}_j}\cos{{\th}_j},\qquad
 \prod^{\ell}_{j=1}dr_j=r^{\ell-1}dr\prod^{\ell-1}_{j=1}\sin^{\ell-1-j}{{\th}_j}d{{\th}_j}.
\end{align*}
 Hence it can be written as 
\[
 g_{\ell}(\mb{R},s)
=(2\pi)^{\ell}\int^{\infty}_{0}\frac{r^{2\ell-1}dr}{(r^2+s^2\bigr)^{\ell+1}}I_{\ell}(r,\mb{R}),
\]
 where $I_{0}(r,\mb{R}):=J_{0}(2\pi r|\mb{R}|)$ and 
\[
 I_{\ell}(r,\mb{R})
:=\underbrace{\int^{\pi/2}_{0}\cdots\int^{\pi/2}_{0}}_{\ell-1}
\prod^{\ell}_{j=1}J_0\biggl(2\pi
 rR_j\bigl(\prod^{j-1}_{k=1}\sin{{\th}_k}\bigr)\cos{{\th}_{j}}\biggr)
\prod^{\ell-1}_{j=1}\sin^{2(\ell-j)-1}{{\th}_j}\cos{{\th}_j}d{{\th}_j}
 \quad (\ell\ge 1).
\]
 By the induction on $\ell$, it is easy to see from the formula
 \eqref{for:Jmn} that 
\[
 I_{\ell}(r,\mb{R})
=\frac{J_{\ell-1}(2\pi r|\mb{R}|)}{(2\pi r|\mb{R}|)^{\ell-1}}.
\]
 Hence, by using \eqref{for:Bessel-int} for $a=s$, $b=2\pi |\mb{R}|$, 
 $\l=\ell-1$ and $\m=\ell$, we have
\[
 g_{\ell}(\mb{R},s)
=\frac{(2\pi)^{\ell}}{(2\pi
 |\mb{R}|)^{\ell-1}}\int^{\infty}_{0}\frac{r^{\ell}J_{\ell-1}(2\pi
 r|\mb{R}|)dr}{(r^2+s^2\bigr)^{\ell+1}}
=\frac{4{\pi}^{\ell+1}|\mb{R}|s^{-1}}{(\ell-1)!}K_{-1}(2\pi|\mb{R}|s).
\]
 Since $K_{-\l}(s)=K_{\l}(s)$, the lemma follows. 
\end{proof}

 By \eqref{for:Fourier-gl}, because $|\mb{y}|=|\mb{R}|$, \Lem{gl} shows
\[
 \Bigl\{\bigl(|\mb{x}+\bsym{\a}|^2+s^2\bigl)^{-\ell-1}\Bigr\}^{\wedge}(\mb{y})
=\frac{2{\pi}^{\ell}}{(\ell-1)!s}e^{-2\pi i\bsym{\a}\mb{y}}\cdot (2{\pi}|\mb{y}|)K_1(2\pi|\mb{y}|s).
\]
 Since $\lim_{s\to 0}sK_{1}(s)=1$, in particular
\[
  \Bigl\{\bigl(|\mb{x}+\bsym{\a}|^2+s^2\bigl)^{-\ell-1}\Bigr\}^{\wedge}(\mb{0})
=\frac{2{\pi}^{\ell}}{(\ell-1)!s^2}.
\]
 The Poisson summation formula therefore yields
\begin{align}
\label{for:PSF-even}
&\ \ \ \frac{\p}{\p s}\Bigl(\frac{1}{2s}\frac{\p}{\p s}\Bigr)^{\ell}\log 
\Bigl(\regprod_{\mb{m}\in \bZ^{2\ell}}
 \bigl(|\mb{m}+\bsym{\a}|^2+s^2\bigr)\Bigr)\nonumber\\ 
&=2s\cdot (-1)^{\ell}{\ell}!\cdot \frac{2{\pi}^{\ell}}{(\ell-1)!s}
\biggl(\frac{1}{s}+\sum_{\mb{n}\in \bZ^{2\ell}\bslO}e^{-2\pi
 i\mb{n}\bsym{\a}}\cdot(2\pi|\mb{n}|)K_1(2\pi|\mb{n}|s)\biggr)\nonumber\\
&=4{\ell}(-1)^{\ell}{\pi}^{\ell}
\biggl(\frac{1}{s}+\sum_{\mb{n}\in \bZ^{2\ell}\bslO}e^{2\pi
 i\mb{n}\bsym{\a}}\cdot(2\pi|\mb{n}|)K_1(2\pi|\mb{n}|s)\biggr).
\end{align}

 In order to integrate this equation, we show the following lemma.

\begin{lem}
 Define a function $Q_{\ell}(s,a)\; (\ell=0,1,2,\ldots)$ of $s$ by  
\[
 Q_{\ell}(s,a):=
\begin{cases}
 \displaystyle{\frac{1}{{\ell}!}s^{2\ell}\log{s}}
  & \quad \textrm{if $a=0$,}\\[10pt]
 (-1)^{\ell+1}(2a^{-1}s)^{\ell}K_{\ell}(as) & \quad \textrm{if $a\ne 0$.}
\end{cases}
\]
 Then, $Q_{\ell}(s,a)$ satisfies the equation  
\begin{align}
\label{for:ladder-Q}
 \frac{\p}{\p s}\Bigl(\frac{1}{2s}\frac{\p}{\p s}\Bigr)^{\ell}
 Q_{\ell}(s,a)
=\begin{cases}
  \displaystyle{\frac{1}{s}} & \quad \textrm{if $a=0$,}\\[10pt]
     aK_1(as) & \quad \textrm{if $a\ne 0$.}
 \end{cases}
\end{align}
\end{lem}
\begin{proof}
 When $a=0$, the assertion is clear. Suppose $a\ne 0$. Since
\[
 K'_{\l}(s)
=\frac{\l}{s}K_{\l}(s)-K_{\l+1}(s)=-\frac{\l}{s}K_{\l}(s)-K_{\l-1}(s),
\]  
 we notice that 
\begin{align}
 \Bigl(\frac{1}{s}\frac{\p}{\p s}\Bigr)^{n}\bigl(s^{\l}K_{\l}(s)\bigr)
&=(-1)^{n}s^{\l-n}K_{\l-n}(s),\label{for:derv-Kn}\\
 K'_{0}(s)
&=-K_1(s).\label{for:derv-K0}
\end{align}
 Hence we have
\begin{align*}
 \frac{\p}{\p s}\Bigl(\frac{1}{2s}\frac{\p}{\p s}\Bigr)^{\ell}Q_{\ell}(s,a)
&=(-1)^{\ell+1}a^{-\ell}\frac{\p}{\p s}
\Bigl(\frac{1}{s}\frac{\p}{\p s}\Bigr)^{\ell}s^{\ell}K_{\ell}(as).
\end{align*}
 Put $u=as$. Since $\frac{\p}{\p s}=a\frac{\p}{\p u}$, we have 
\begin{align*}
 \frac{\p}{\p s}\Bigl(\frac{1}{2s}
\frac{\p}{\p s}\Bigr)^{\ell}Q_{\ell}(s,a)
&=(-1)^{\ell+1}a^{-\ell+1}\frac{\p}{\p u}
 \bigl(a^2\frac{1}{u}\frac{\p}{\p u}\bigr)^{\ell}(a^{-1}u)^{\ell}K_{\ell}(u)\\
&=(-1)^{\ell+1}a\frac{\p}{\p u}
 \bigl(\frac{1}{u}\frac{\p}{\p u}\bigr)^{\ell}u^{\ell}K_{\ell}(u)\\
&=(-1)^{\ell+1}a\frac{\p}{\p u}(-1)^{\ell}K_0(u)\\
&=aK_1(as).
\end{align*}
 Here we have used the equations \eqref{for:derv-Kn} and
 \eqref{for:derv-K0}.
\end{proof}

 By this lemma, it follows immediately from \eqref{for:PSF-even} that 
\begin{align*} 
\frac{\p}{\p s}
\Bigl(\frac{1}{2s}\frac{\p}{\p s}\Bigr)^{\ell}
& \Bigg[\log 
\Bigl(\regprod_{\mb{m}\in \bZ^{2\ell}}
 \bigl(|\mb{m}+\bsym{\a}|^2+s^2\bigr)\Bigr)\\ 
&-\biggl\{\frac{4(-1)^{\ell}{\pi}^{\ell}}{(\ell-1)!}s^{2\ell}\log{s}-4{\ell}s^{\ell}\sum_{\mb{n}\in
 \bZ^{2\ell}\bslO}\frac{1}{|\mb{n}|^{\ell}}e^{2\pi i\mb{n}\bsym{\a}}K_{\ell}(2\pi|\mb{n}|s)\biggr\}\Bigg]=0.
\end{align*}
\end{enumerate}
 This completes the proof of the theorem.
\end{proof} 

 Using $\sum_{\mb{n}\in \bZ^{\n}\bslO}=\sum^{\infty}_{n=1}\sum_{\mb{m}\in\bZ^{\n}, |\mb{m}|^2=n}$,
 we immediately obtain from \Thm{DetLaplacian} the following

\begin{cor}
\label{cor:series}
 Put 
\begin{align*}
 r_{\n}(n,\bsym{\a})
:=\sum_{\mb{m}\in \bZ^{\n},\ |\mb{m}|^2=n}e^{2\pi i\mb{m}\bsym{\a}}.
\end{align*}
 For $\Re(s)>0$, we have
\begin{align}
 \quad\det(\D_{2\ell+1,\bsym{\a}}+s^2)
&\cong_{2\ell}\exp\biggl(-\frac{(-2\pi)^{\ell+1}}{(2\ell+1)!!}s^{2\ell+1}
-\sum^{\infty}_{n=1}\sum^{\ell}_{k=0}\frac{r_{2\ell+1}(n,\bsym{\a})c^{(\ell)}_k}{n^{\frac{\ell+1}{2}}(2\pi\sqrt{n})^k}s^{\ell-k}e^{-2\pi\sqrt{n}s}\biggr),\\
 \quad\det(\D_{2\ell,\bsym{\a}}+s^2)
&\cong_{2\ell}\exp\biggl(\frac{4(-1)^{\ell}{\pi}^{\ell}}{(\ell-1)!}s^{2\ell}\log{s}-4{\ell}s^{\ell}
\sum^{\infty}_{n=1}\frac{r_{2\ell}(n,\bsym{\a})}{n^{\frac{\ell}{2}}}K_{\ell}(2\pi\sqrt{n}s)\biggr).  
\end{align} \qed
\end{cor}

\begin{remark}
 The series $r_{\n}(n,\bsym{\a})$ have been studied, for instance in
 \cite{BleherBourgain1996}, \cite{BleherDyson1994}, etc, in connection
 with the fluctuations of the number of lattice points inside a large
 sphere centered at $\bsym{\a}\in\bR^{\n}$. In particular, quite
 recently, the limit distribution of the mean value of the square sum
 $X^{-\frac{\n}{2}}\sum_{n \le X}|r_{\n}(n,\bsym{\a})|^2$ for $X\to \infty$
 has been obtained in \cite{Marklof2005} explicitly when $\bsym{\a}$
 satisfies a certain diophantine condition.
\end{remark}

\section{Arithmetic functions arising from $\bsym{\log L_{\pp}(s,\bsym{\a})}$}  

 By the definition of $L_{\pp}(s,\bsym{\a};\n)$ we can calculate 
\[
 \log L_{\pp}(s,\bsym{\a};\n)
=\sum_{\gcd\mb{p}=1}\sum^{\infty}_{\ell=1}
\frac{e^{2\pi i\ell\mb{p}\bsym{\a}}e^{-s\ell|\mb{p}|}}{\ell}
=\sum^{\infty}_{n=1}\sum_{|\mb{m}|^2=n}(\gcd\mb{m})^{-1}
e^{2\pi i\mb{m}\bsym{\a}}e^{-s\sqrt{n}}.   
\]
 Hence, if we put
\begin{equation}
 M_{\n}(n,\bsym{\a},x)
=\sum_{\mb{m}\in\bZ^{\n},\ |\mb{m}|^2=n}(\gcd\mb{m})^{-x}e^{2\pi i\mb{m}\bsym{\a}}\end{equation}
 we have 
\begin{align}
\label{for:logLpp}
 \log L_{\pp}(s,\bsym{\a};\n)
=\sum^{\infty}_{n=1} M_{\n}(n,\bsym{\a},1)e^{-s\sqrt{n}}.
\end{align}
 We remark that $M_{\n}(n,\bsym{\a},0)=r_{\n}(n,\bsym{\a})$, where
 $r_{\n}(n,\bsym{\a})$ is defined in Corollary\,\ref{cor:series}. In
 \cite{TorusSelbergZeta} we prove that $M_{\n}(n,\mb{0},x)$ is a 
 multiplicative function (w.r.t. the variable $n$) when $\n=2$ and study
 an asymptotic distribution of the average. The aim of this section is
 to generalize the results in \cite{TorusSelbergZeta} to the cases of
 $\n=4,6$ and $8$. 

 We first notice that  
\begin{equation}
\label{for:Mn}
 M_{\n}(n,\bsym{\a},x)
=\sum^{\infty}_{\ell=1}
 \ell^{-x}\sum_{\gcd\mb{m}=\ell,|\mb{m}|^2=n}e^{2\pi i\mb{m}\bsym{\a}} 
=\sum_{{\ell}^2|n}g_x({\ell}^2)\w{r}_{\n}\bigl(\frac{n}{\ell^2},\ell\bsym{\a}\bigr),\end{equation}
 where 
\begin{align*}
 g_x(n)
:=
\begin{cases}
 n^{-\frac{x}2} & \quad \text{if $n$ is square},\\
              0 & \quad \text{otherwise},
\end{cases}
\qquad
\textrm{and}
\qquad
 \w{r}_{\n}(n,\bsym{\a})
:=\sum_{\mb{m}\in\bZ^{\n},\ |\mb{m}|^2=n \atop \gcd\mb{m}=1}
e^{2\pi i\mb{m}\bsym{\a}}.
\end{align*}
 We note that the function $g_x(n)$ is multiplicative. Hence if
 $\bsym{\a}=\mb{0}$, we have $M_{\n}(n,x)=g_x * \w{r}_{\n}(n)$ where
 $M_{\n}(n,x):=M_{\n}(n,\mb{0},x)$ and
 $\w{r}_{\n}(n):=\w{r}_{\n}(n,\mb{0})$. Notice also that   
\begin{equation}
\label{for:Rr}
 r_{\n}(n,\bsym{\a})
=\sum_{{\ell}^2|n}\w{r}_{\n}\bigl(\frac{n}{{\ell}^2},\ell\bsym{\a}\bigr),
\end{equation}

 Let 
\[
 D_{\n}(s;\bsym{\a},x):=\sum^{\infty}_{n=1}M_{\n}(n,\bsym{\a},x)n^{-s}.
\]
 From \eqref{for:Mn}, we have
\begin{align}
\label{for:DDD}
 D_{\n}(s;\bsym{\a},x)
&=\sum^{\infty}_{n=1}\sum_{{\ell}^2|n}g_x({\ell}^2)\w{r}_{\n}\bigl(\frac{n}{\ell^2},\ell\bsym{\a}\bigr)n^{-s}=\sum^{\infty}_{\ell=1}\sum^{\infty}_{k=1}g_x({\ell}^2)\w{r}_{\n}\bigl(\frac{k{\ell}^2}{\ell^2},\ell\bsym{\a}\bigr)(k{\ell}^2)^{-s}\nonumber\\
&=\sum^{\infty}_{\ell=1}g_x({\ell}^2){\ell}^{-2s}\sum^{\infty}_{k=1}\w{r}_{\n}(k,\ell\bsym{\a}\bigr)k^{-s}. 
\end{align}
 We further define the functions
\[
 \cL_{\n}(s;\bsym{\a}):=\sum^{\infty}_{n=1}r_{\n}(n,\bsym{\a})n^{-s}
 \quad \textrm{and} \quad  
 \w{\cL}_{\n}(s;\bsym{\a}):=\sum^{\infty}_{n=1}\w{r}_{\n}(n,\bsym{\a})n^{-s}
\]
 (for the study of the function $\cL_{\n}(s;\mb{0})$, see
 \cite{BorweinChoi2003}). Then we have from \eqref{for:Rr} and
 \eqref{for:DDD} that
\[
\cL_{\n}(s;\bsym{\a})
=\sum^{\infty}_{\ell=1}{\ell}^{-2s}\w{\cL}_{\n}(s;\ell\bsym{\a})
\qquad \textrm{and} \qquad 
  D_{\n}(s;\bsym{\a},x)
=\sum^{\infty}_{\ell=1}g_x({\ell}^2){\ell}^{-2s}\w{\cL}_{\n}(s;\ell\bsym{\a}). 
\]
 In particular, we put $\cL_{\n}(s):=\cL_{\n}(s;\mb{0})$ and
 $\w{\cL}_{\n}(s):=\w{\cL}_{\n}(s;\mb{0})$ when $\bsym{\a}=\mb{0}$.
 Then we have  
\begin{equation}
\label{for:a0}
 \cL_{\n}(s)=\z(2s)\w{\cL}_{\n}(s)
\qquad \textrm{and} \qquad 
  D_{\n}(s;x)=\z(x+2s)\w{\cL}_{\n}(s),
\end{equation}
 because
\[
 \sum^{\infty}_{\ell=1}g_x({\ell}^2){\ell}^{-2s}
=\sum^{\infty}_{\ell=1}({\ell}^2)^{-\frac{x}{2}}{\ell}^{-2s}
=\sum^{\infty}_{\ell=1}{\ell}^{-x-2s}
=\z(x+2s).
\]
 From \eqref{for:a0}, we obtain
\begin{equation}
\label{for:D}
 D_{\n}(s;x)=\z(x+2s)\z(2s)^{-1}\cL_{\n}(s).
\end{equation}
 Note that $r_{\n}(n,\mb{0})=r_{\n}(n)$. Thus, a similar discussion
 performed in \cite{TorusSelbergZeta} gives the following asymptotic
 average of $M_{\n}(n,\mb{0},x)$ for $\n=2,4,6$ and $8$.

\begin{example}
[The case $\bsym{\n=2}$ (studied in \cite{TorusSelbergZeta})]
 From \eqref{for:r2}, we have 
\begin{align*}
 \cL_{2}(s)&=4\z(s) L_{-4}(s). 
\end{align*}
 where $L_{-4}(s):=\sum^{\infty}_{n=1}\c_{-4}(n)n^{-s}$. Hence it
 follows from \eqref{for:D} that 
\[
 D_{2}(s;x)=4\z(x+2s)\z(2s)^{-1}\z(s)L_{-4}(s).
\] 
 Therefore, by the Tauberian theorem (see \cite{MurtyMurty1997}), we
 obtain  
\[
 \sum_{n\le X}M_2(n,x)=
\begin{cases}
 \displaystyle{\biggl(\frac{2\z\bigl(\frac{-x+1}{2}\bigr)L_{-4}\bigl(\frac{-x+1}{2}\bigr)}{\z(-x+1)}+o(1)}\biggr)X^{\frac{-x+1}{2}}
 & \quad \textrm{if} \quad x<-1,\\[10pt]
 \displaystyle{\Bigl(\frac{3}{\pi}+o(1)\Bigr)X\log{X}}
 & \quad \textrm{if} \quad x=-1,\\[10pt]
 \displaystyle{\Bigl(\frac{2\cdot 3\z(x+2)}{\pi}+o(1)\Bigr)X}
 & \quad \textrm{if} \quad x>-1
\end{cases}
\]
 as $X\to \infty$, since $\z(2)={\pi}^2/(2\cdot 3)$ and
 $L_{-4}(1)={\pi}/2^2$. In particular, we have
\[
 \sum_{n\le X}M_2(n,1)
=\Bigl(\frac{2\cdot 3}{\pi}\z(3)+o(1)\Bigr)X \qquad (X\to \infty).
\] 
\end{example}

\begin{remark}
 The value $6\z(2+x)/{\pi}^2$ for the case of $x>-1$ in Corollary\,$4.4$
 in \cite{TorusSelbergZeta} is incorrect and should be $6\z(2+x)/\pi$. 
\end{remark}

\begin{remark}
 From the example above, we have
\[
 \lim_{X\to \infty}\frac{1}{X}\sum_{n \le X}M_2(n,1)
=\pi\cdot \frac{\z(3)}{\z(2)}.
\]
 The appearance of the ratio of $\z(3)/\z(2)$ can be also found at the
 study of the mean square limit for lattice points in the $3$-dim sphere
 in \cite{Jarnik1940} and \cite{BleherDyson1994-2}.
\end{remark}

\begin{example}
[The case $\bsym{\n=4}$]
 One can calculate from \eqref{for:r4} that
\begin{align*}
 \cL_{4}(s)&=8(1-4^{1-s})\z(s)\z(s-1)
\end{align*}
 and hence
\[
 D_{4}(s;x)=8(1-4^{1-s})\z(x+2s)\z(2s)^{-1}\z(s)\z(s-1).
\] 
 Therefore it holds that 
\[
 \sum_{n\le X}M_4(n,x)=
\begin{cases}
 \displaystyle{\biggl(\frac{4(1-4^{\frac{x+1}{2}})\z\bigl(\frac{-x+1}{2}\bigr)\z\bigl(\frac{-x-1}{2}\bigr)}{\z(-x+1)}+o(1)}\biggr)X^{\frac{-x+1}{2}}
 & \quad \textrm{if} \quad x<-3,\\[10pt]
 \displaystyle{\Bigl(\frac{3^2 \cdot 5}{{\pi}^2}+o(1)\Bigr)X^2\log{X}}
 & \quad \textrm{if} \quad x=-3,\\[10pt]
 \displaystyle{\Bigl(\frac{2\cdot 3^2\cdot 5\z(x+4)}{{\pi}^2}+o(1)\Bigr)X^2}
 & \quad \textrm{if} \quad x>-3
\end{cases}
\]
 as $X\to \infty$, since $\z(4)={\pi}^4/(2\cdot 3^2\cdot 5)$. In
 particular, we have   
\[
 \sum_{n\le X}M_4(n,1)
=\Bigl(\frac{2\cdot 3^2\cdot 5}{{\pi}^2}\z(5)+o(1)\Bigr)X^2 \qquad (X\to\infty).\]
\end{example}

\begin{example}
[The case $\bsym{\n=6}$]
 
 It is known (see, e.g., \cite{BorweinChoi2003}) that
\begin{align*}
 r_{6}(n)
&=16\sum_{m|n}\c_{-4}\Bigl(\frac{n}{m}\Bigr)m^2-4\sum_{m|n}\c_{-4}(m)m^2,\\
 \cL_{6}(s)
&=16\z(s-2)L_{-4}(s)-4\z(s)L_{-4}(s-2).
\end{align*}
 Hence we have
\[
 D_{6}(s;x)=\z(x+2s)\z(2s)^{-1}\bigl(16\z(s-2)L_{-4}(s)-4\z(s)L_{-4}(s-2)\bigr).\]
 Therefore we obtain 
\[
 \sum_{n\le X}M_6(n,x)=
\begin{cases}
 \displaystyle{\biggl(\frac{8\z\bigl(\frac{-x-3}{2}\bigr)L_{-4}\bigl(\frac{-x+1}{2}\bigr)-2\z\bigl(\frac{-x+1}{2}\bigr)L_{-4}\bigl(\frac{-x-3}{2}\bigr)}{\z(-x+1)}+o(1)}\biggr)X^{\frac{-x+1}{2}}
 & \quad \textrm{if} \quad x<-5,\\[10pt]
 \displaystyle{\Bigl(\frac{3^3\cdot 5\cdot 7}{2^2{\pi}^3}+o(1)\Bigr)X^3\log{X}}
 & \quad \textrm{if} \quad x=-5,\\[10pt]
 \displaystyle{\Bigl(\frac{3^3\cdot 5\cdot 7\z(x+6)}{2{\pi}^3}+o(1)\Bigr)X^3}
 & \quad \textrm{if} \quad x>-5
\end{cases}
\]
 as $X\to \infty$, since $\z(6)={\pi}^6/(3^3\cdot 5\cdot 7)$ and
 $L_{-4}(3)={{\pi}^3}/2^5$. In particular, we have   
\[
 \sum_{n\le X}M_6(n,1)
=\Bigl(\frac{3^3\cdot 5\cdot 7}{2{\pi}^3}\z(7)+o(1)\Bigr)X^3 \qquad (X\to \infty).
\]
\end{example}

\begin{example}
[The case $\bsym{\n=8}$]
 We have from \eqref{for:r8} that 
\begin{align*}
 \cL_{8}(s)&=16(1-2^{1-s}+4^{2-s})\z(s)\z(s-3).
\end{align*}
 Hence
\[
 D_{8}(s;x)=16(1-2^{1-s}+4^{2-s})\z(x+2s)\z(2s)^{-1}\z(s)\z(s-3).
\] 
 Therefore we obtain 
\[
 \sum_{n\le X}M_8(n,x)=
\begin{cases}
 \displaystyle{\biggl(\frac{8(1-2^{\frac{x+1}{2}}+4^{\frac{x+3}{2}})\z\bigl(\frac{-x+1}{2}\bigr)\z\bigl(\frac{-x-5}{2}\bigr)}{\z(-x+1)}+o(1)}\biggr)X^{\frac{-x+1}{2}}
 & \quad \textrm{if} \quad x<-7,\\[10pt]
 \displaystyle{\Bigl(\frac{3^2 \cdot 5^2\cdot 7}{2{\pi}^4}+o(1)\Bigr)X^4\log{X}} & \quad \textrm{if} \quad x=-7,\\[10pt]
 \displaystyle{\Bigl(\frac{3^2 \cdot 5^2\cdot 7\z(x+8)}{{\pi}^4}+o(1)\Bigr)X^4}
 & \quad \textrm{if} \quad x>-7
\end{cases}
\]
 as $X\to \infty$, since 
 $\z(8)={\pi}^8/(2\cdot 3^3\cdot 5^2\cdot 7)$. In particular, we have   
\[
 \sum_{n\le X}M_8(n,1)
=\Bigl(\frac{3^2 \cdot 5^2\cdot 7}{{\pi}^4}\z(9)+o(1)\Bigr)X^4 \qquad
       (X\to\infty).
\]
\end{example}

 From the observations above, we strongly expect the following statement
 can be true.

\begin{conj}
\label{conj:arith}
 For any $\ell\in\bN$, there exists $\b_{\ell}\in\bQ$ such that 
\[
 \sum_{n\le X}M_{2\ell}(n,1)
=\Bigl(\frac{\b_{\ell}}{{\pi}^{\ell}}\z(2\ell+1)+o(1)\Bigr)X^{\ell} \qquad
       (X\to\infty).
\]
\end{conj}

\bigbreak 

\noindent
 {\it Additional reports (September 2007).}\ We have obtained the
 conjecture above affirmatively. In fact, for any $\n\in\bN$, it holds
 that 
\begin{equation}
\label{for:any-nu}
 \sum_{n\le X}M_{\n}(n,x)=
\begin{cases}
 \displaystyle{\biggl(\frac{\cL_{\n}\bigl(\frac{1-x}{2}\bigr)}{\z(-x+1)}+o(1)}\biggr)X^{\frac{-x+1}{2}}
 & \quad \textrm{if} \quad x<1-\n,\\[10pt]
 \displaystyle{\Bigl(\frac{\pi^{\frac{\n}{2}}}{2\z(\n)\G(\frac{\n}{2})}+o(1)\Bigr)X^{\frac{\n}{2}}}\log{X} & \quad \textrm{if} \quad x=1-\n,\\[10pt]
 \displaystyle{\Bigl(\frac{\pi^{\frac{\n}{2}}}{\z(\n)\G(\frac{\n}{2})}\z(\n+x)+o(1)\Bigr)X^{\frac{\n}{2}}}
 & \quad \textrm{if} \quad x>1-\n
\end{cases}
\end{equation}  
 as $X\to\infty$. In particular, for $\ell\in\bN$, we have
\begin{align}
\label{for:even}
 \sum_{n\le X}M_{2\ell}(n,1)
&=\Bigl(\frac{(-1)^{\ell+1}(2\ell)!}{(\ell-1)!2^{2\ell-1}B_{2\ell}}\frac{\z(2\ell+1)}{\pi^{\ell}}+o(1)\Bigr)X^{\ell}
 \qquad (X\to\infty),\\
\label{for:odd}
 \sum_{n\le X}M_{2\ell+1}(n,1)
&=\Bigl(\frac{(-1)^{\ell}2^{3\ell+1}B_{2\ell+2}}{(2\ell-1)!!(2\ell+2)!}\frac{\pi^{3\ell+2}}{\z(2\ell+1)}+o(1)\Bigr)X^{\ell+\frac{1}{2}}
\qquad (X\to\infty),
\end{align}
 where $B_{n}$ is the Bernoulli number (namely, the number $\b_{\ell}$
 in Conjecture~\ref{conj:arith} is given by 
 $\b_{\ell}=\frac{(-1)^{\ell+1}(2\ell)!}{(\ell-1)!2^{2\ell-1}B_{2\ell}}\in\bQ$). Actually, the formula \eqref{for:any-nu} follows from the application
 of the Tauberian theorem to \eqref{for:D} with
 $\cL_{\n}(s)=Z(s,1_{\n})$ where $1_{\n}$ is the identity matrix of
 order $\n$, $Z(s,A):=\sum_{\bsym{m}\in\bZ^{\n}}A[\bsym{m}]^{-s}$ is the
 Epstein zeta function attached to the positive definite symmetric
 matrix $A$ and $A[\bsym{x}]:={}^t\bsym{x}A\bsym{x}$ for
 $\bsym{x}\in\bR^{\n}$. Note that $Z(s,A)$ converges absolutely for
 $\Re(s)>\n/2$ and admit a meromorphic continuation to the whole plane
 $\bC$ with a simple pole at $s=\n/2$ with residue
 $\pi^{\frac{\n}{2}}/(\sqrt{\det{A}}\cdot\G(\frac{\n}{2}))$ (see 
 \cite{Terras1985}). Moreover, using
 $\z(2n)=(-1)^{n+1}2^{2n-1}B_{2n}/(2n)!$ for $n\in\bN$, one can easily
 obtain the formulas \eqref{for:even} and \eqref{for:odd} from
 \eqref{for:any-nu}.


\bigskip

\noindent
\textsc{Nobushige KUROKAWA}\\
 Department of Mathematics, Tokyo Institute of Technology.\\
 Oh-okayama Meguro, Tokyo, 152-0033 JAPAN.\\
 \texttt{kurokawa@math.titech.ac.jp}\\

\noindent
\textsc{Masato WAKAYAMA}\\
 Faculty of Mathematics, Kyushu University.\\
 Hakozaki, Fukuoka, 812-8581 JAPAN.\\
 \texttt{wakayama@math.kyushu-u.ac.jp}\\

\noindent
 \textsc{Yoshinori YAMASAKI}\\
 Graduate School of Mathematics, Kyushu University.\\
 Hakozaki, Fukuoka, 812-8581 JAPAN.\\
 \texttt{yamasaki@math.kyushu-u.ac.jp}\\

\end{document}